\documentclass[10pt]{article}
\usepackage{latexsym}
\usepackage{amsfonts}
\usepackage{enumerate}
\usepackage{multicol}
\usepackage{graphicx}
\usepackage{amssymb}
\usepackage{amsmath}
\usepackage{epic}

\title{On the Minimum Size of Signed Sumsets in Elementary Abelian Groups \\[.4in]}

\author{B\'{e}la Bajnok\footnote{Corresponding author} \\[.1in] {\small Department of Mathematics, Gettysburg College} \\
{\small 300 N. Washington Street, Gettysburg, PA 17325-1486 USA} \\{\small E-mail:  bbajnok@gettysburg.edu} \\ [.2in]
and \\[.2in]
Ryan Matzke \\[.1in] {\small Department of Mathematics, Gettysburg College} \\
{\small 300 N. Washington Street, Gettysburg, PA 17325-1486 USA} \\{\small E-mail:  matzry01@gettysburg.edu} \\ [.2in]
 \\[.4in]}

\date{December 2, 2014}

\newtheorem{thm}{Theorem}

\newtheorem{cor}[thm]{Corollary}
\newtheorem{prop}[thm]{Proposition}
\newtheorem{conj}[thm]{Conjecture}

\begin{document}

\maketitle

\begin{abstract}

For a finite abelian group $G$ and positive integers $m$ and $h$, we let 
$$\rho(G, m, h) = \min \{ |hA| \; : \; A \subseteq G, |A|=m\}$$ and
$$\rho_{\pm} (G, m, h) = \min \{ |h_{\pm} A| \; : \; A \subseteq G, |A|=m\},$$ where $hA$ and $h_{\pm} A$ denote the $h$-fold sumset and the $h$-fold signed sumset of $A$, respectively.  The study of $\rho(G, m, h)$ has a 200-year-old history and is now known for all $G$, $m$, and $h$.  In previous work we provided an upper bound for $\rho_{\pm} (G, m, h)$ that we believe is exact, and proved that $\rho_{\pm} (G, m, h)$ agrees with $\rho (G, m, h)$ when $G$ is cyclic.  Here we study $\rho_{\pm} (G, m, h)$ for elementary abelian groups $G$; in particular, we determine all values of $m$ for which $\rho_{\pm} (\mathbb{Z}_p^2, m, 2)$ equals $\rho (\mathbb{Z}_p^2, m, 2)$ for a given prime $p$.

\end{abstract}



\thispagestyle{empty}

\section{Introduction}

Let $G$ be a finite abelian group written with additive notation, let $m$ be a positive integer with $m \leq |G|$, and let $h$ be a nonnegative integer.  In \cite{BajMat:2014a}, we introduced the function
$$\rho_{\pm}(G, m, h) = \min \{ |h_{\pm}A|  \; : \; A \subseteq G, |A|=m\},$$ where
$$h_{\pm} A=\{ \Sigma_{i=1}^m \lambda_i a_i \; : \; (\lambda_1,\dots,\lambda_m) \in \mathbb{Z}^m, \; \Sigma_{i=1}^m |\lambda_i|=h\}$$ is the 
$h$-fold {\em signed sumset} of an $m$-subset $A=\{a_1, \dots, a_m\}$ of $G$
 (as usual, $|S|$ denotes the size of the finite set $S$).  The function $\rho_{\pm}(G, m, h)$ is the analogue of the well-known 
$$\rho(G, m, h) = \min \{ |hA| \; : \; A \subseteq G, |A|=m\},$$ where  
$$hA=\{ \Sigma_{i=1}^m \lambda_i a_i \; : \; (\lambda_1,\dots,\lambda_m) \in \mathbb{N}_0^m, \; \Sigma_{i=1}^m \lambda_i=h\}$$ is the usual $h$-fold {\em sumset} of $A$.  

Signed sumsets have already been studied in the past: For example, in \cite{BajRuz:2003a}, the first author and Ruzsa investigated the {\em independence number} of a subset $A$ of $G$, defined as the maximum value of $t \in \mathbb{N}$ for which $$0 \not \in \cup_{h=1}^t h_{\pm}A$$ (see also \cite{Baj:2000a} and \cite{Baj:2004a}); and in \cite{KloLev:2003a}, Klopsch and Lev discussed the {\em diameter} of $G$ with respect to $A$, defined as the minimum value of $s \in \mathbb{N}$ for which $$\cup_{h=0}^s h_{\pm}A=G$$ (see also \cite{KloLev:2009a}).  The independence number of $A$ in $G$ quantifies the ``degree'' to which $A$ is linearly independent in $G$, while the diameter of $G$ with respect to $A$ measures how ``effectively'' $A$ generates $G$ (if at all).  While research on minimum sumset size goes back to the work of Cauchy and is now known for all $G$, $m$, and $h$, to the best of our knowledge, \cite{BajMat:2014a} is the first systematic study of the minimum size of signed sumsets.  In this paper we continue our work and consider $\rho_{\pm}(G, m, h)$ for elementary abelian groups $G$.

Let us review what we need to know about $\rho(G, m, h)$.  It has been over two hundred years since Cauchy \cite{Cau:1813a} found the minimum possible size of
$$A+B=\{a+b \; : \; a \in A,\; b \in B \}$$ among subsets $A$ and $B$ of the cyclic group $\mathbb{Z}_p$ of given sizes.  (Here and elsewhere in the paper $p$ denotes a positive prime.)  Over a hundred years later, Davenport \cite{Dav:1935a} (cf.~\cite{Dav:1947a}) rediscovered Cauchy's result, which is now known as the Cauchy--Davenport Theorem: 

\begin{thm}[Cauchy--Davenport Theorem] \label{Cauchy--Davenport}
If $A$ and $B$ are nonempty subsets of the group $\mathbb{Z}_p$ of prime order $p$, then 
$$|A+B| \geq \min \{p, |A|+|B|-1\}.$$
\end{thm}
It can easily be seen that the bound is tight for all values of $|A|$ and $|B|$, and thus
$$ \rho (\mathbb{Z}_p, m, 2)=\min\{p,2m-1\}.$$

Relatively recently, $ \rho (G, m, h)$ was finally evaluated for all parameters by Plagne \cite{Pla:2006a} (see also \cite{Pla:2003a}, \cite{EliKer:2007a}, and  \cite{EliKerPla:2003a}) in 2003.  To state the result, we introduce the function
$$u(n,m,h)=\min \{f_d (m,h) \; : \; d \in D(n)\},$$ where $n$, $m$, and $h$ are positive integers, $D(n)$ is the set of positive divisors of $n$, and
$$f_d(m,h)=\left(h\left \lceil m/h \right \rceil-h +1 \right) \cdot d.$$  
(Here $u(n,m,h)$ is a relative of the Hopf--Stiefel function used also in topology and bilinear algebra; see, for example, \cite{EliKer:2005a}, \cite{Kar:2006a}, \cite{Pla:2003a}, and  \cite{Sha:1984a}.)

\begin{thm} [Plagne; cf.~\cite{Pla:2006a}] \label{value of u}
Let $n$, $m$, and $h$ be positive integers with $m \leq n$.  For any abelian group $G$ of order $n$ we have
$$\rho (G, m, h)=u(n,m,h).$$
\end{thm}  

Let us turn now to $\rho_{\pm} (G, m, h)$.  It is easy to see that $\rho_{\pm} (G,1,h)$ and $\rho_{\pm} (G,m,0)$ both equal $1$ and that $\rho_{\pm} (G,m,1)$ equals $m$ for all $G$, $m$, and $h$.  (To see the last equality, it suffices to verify that one can always find a {\em symmetric} subset of size $m$ in $G$, that is, an $m$-subset $A$ of $G$ for which $A=-A$.)   Therefore, from now on, we assume that $m \geq 2$ and $h \geq 2$.  

Perhaps surprisingly, we find that, while the $h$-fold signed sumset of a given set is generally much larger than its sumset, $\rho_{\pm} (G, m, h)$ often agrees with $\rho (G, m, h)$; in particular, this is always the case when $G$ is cyclic:  

\begin{thm} [Cf.~\cite{BajMat:2014a}] \label{cyclic} For all positive integers $n$, $m$, and $h$, we have 
$$\rho_{\pm} (\mathbb{Z}_n, m, h)= \rho (\mathbb{Z}_n, m, h).$$
\end{thm}

The situation seems considerably more complicated for noncyclic groups: in contrast to $\rho (G, m,h)$, the value of $\rho_{\pm} (G, m,h)$ depends on the structure of $G$ rather than just the order $n$ of $G$.  

Observe that by Theorem \ref{value of u}, we have the lower bound
$$\rho_{\pm} (G, m,h) \geq u(n,m,h)=\min \{f_d (m,h) \; : \; d \in D(n)\}.$$ In \cite{BajMat:2014a}, we proved that with a certain subset $D(G,m)$ of $D(n)$, we have
$$\rho_{\pm} (G, m,h) \leq u_{\pm} (G,m,h)=\min \{f_d (m,h) \; : \; d \in D(G,m)\};$$ here $D(G,m)$ is defined in terms of the {\em type} $(n_1,\dots,n_r)$ of $G$, that is, via  integers $n_1,\dots,n_r$ such that $n_1 \geq 2$, $n_i$ divides $n_{i+1}$ for each $i \in \{1,\dots, r-1\}$, and for which $G$ is isomorphic to  the invariant product  
$$\mathbb{Z}_{n_1} \times \cdots \times \mathbb{Z}_{n_r}.$$
Namely, we proved the following result:

\begin{thm} [Cf.~\cite{BajMat:2014a}]  \label{u pm with f}  The minimum size of the $h$-fold signed sumset of an $m$-subset of a group $G$ of type $(n_1,\dots,n_r)$ satisfies
$$\rho_{\pm} (G, m,h) \leq  u_{\pm} (G,m,h),$$
where $$u_{\pm} (G,m,h)=\min \{f_d (m,h) \; : \; d \in D(G,m) \}$$ with $$D(G,m)=\{d \in D(n) \; : \; d= d_1 \cdots d_r, d_1 \in D(n_1), \dots, d_r \in D(n_r), dn_r \geq d_rm \}.$$
\end{thm}
Observe that, for cyclic groups of order $n$, $D(G,m)$ is simply $D(n)$. 

Additionally, we believe that $u_{\pm} (G,m,h)$ actually yields the exact value of $\rho_{\pm} (G,m,h)$ in all cases except for one very special situation (which occurs only when $h=2$).  In particular, we made the following conjecture:

\begin{conj} [Cf.~\cite{BajMat:2014a}] \label{conj for rho pm}
Suppose that $G$ is an abelian  group of order $n$ and type $(n_1, \dots, n_r)$.

If $h \geq 3$, then $$\rho_{\pm} \left(G, m, h   \right) =u_{\pm}(G,m,h).$$

If each odd divisor of $n$ is less than $2m$, then $$\rho_{\pm} \left(G, m, 2  \right) =u_{\pm}(G,m,2).$$

If there are odd divisors of $n$ greater than $2m$, let $d_m$ be the smallest one.  We then have
$$\rho_{\pm} \left(G, m, 2   \right) = \min\{u_{\pm}(G,m,2), d_m-1\}.$$
\end{conj}

We will need to use the following ``inverse type'' result from \cite{BajMat:2014a} regarding subsets that achieve $\rho_{\pm} \left(G, m, h   \right)$.  Given a group $G$ and a positive integer $m \leq |G|$, we define a certain collection ${\cal A}(G,m)$ of $m$-subsets of $G$.  
We let 
\begin{itemize}
  \item $\mathrm{Sym}(G,m)$ be the collection of {\em symmetric} $m$-subsets of $G$, that is, $m$-subsets $A$ of $G$ for which $A=-A$;
  \item $\mathrm{Nsym}(G,m)$ be the collection of {\em near-symmetric} $m$-subsets of $G$, that is, $m$-subsets $A$ of $G$ that are not symmetric, but for which $A\setminus \{a\}$ is symmetric for some $a \in A$; 
  \item $\mathrm{Asym}(G,m)$ be the collection of {\em asymmetric} $m$-subsets of $G$, that is, $m$-subsets $A$ of $G$ for which $A \cap (-A)=\emptyset$.  
  \end{itemize}  We then let 
$${\cal A}(G,m)=\mathrm{Sym}(G,m) \cup \mathrm{Nsym}(G,m)\cup \mathrm{Asym}(G,m).$$  In other words, ${\cal A}(G,m)$ consists of those $m$-subsets of $G$ that have exactly $m$, $m-1$, or $0$ elements whose inverse is also in the set. 

\begin{thm} [Cf.~\cite{BajMat:2014a}] \label{symmetry thm}
For every $G$, $m$, and $h$, we have
$$\rho_{\pm} (G,m,h)= \min \{|h_{\pm} A| \; : \; A \in {\cal A}(G,m)\}.$$

\end{thm}
We should add that each of the three types of sets are essential as can be seen by examples (cf.~\cite{BajMat:2014a}).   

Our goal in this paper is to investigate $\rho_{\pm} (G,m,h)$ for elementary abelian groups $G$.  In particular, we wish to classify all cases for which
$$\rho_{\pm} (\mathbb{Z}_p^r, m, h) = \rho (\mathbb{Z}_p^r, m, h),$$
where $p$ denotes a positive prime and $r$ is a positive integer.  By Theorem \ref{cyclic}, we assume that $r \geq 2$, and, since obviously   
$$\rho_{\pm} (\mathbb{Z}_2^r, m, h) = \rho (\mathbb{Z}_2^r, m, h)$$ for all $m$, $h$, and $r$, we will also assume that $p \geq 3$.  

Let us first exhibit a sufficient condition for $\rho_{\pm} (\mathbb{Z}_p^r, m, h)$ to equal $\rho(\mathbb{Z}_p^r, m, h)$.  When $p \leq h$, our result is easy to state; we will prove the following:

\begin{thm} \label{p leq h}

If $p \leq h$, then for all values of $1 \leq m \leq p^r$ we have $$\rho_{\pm} (\mathbb{Z}_p^r, m, h) = \rho (\mathbb{Z}_p^r, m, h).$$

\end{thm}

The case $h \leq p-1$ is more complicated and delicate.  In order to state our results, we will need to introduce some notations.  Suppose that $m \geq 2$ is a given positive integer.  First, we let $k$ be the maximal integer for which
$$p^k  +\delta \leq hm-h+1,$$  where 
$\delta=0$ if $p-1$ is divisible by $h$, and $\delta=1$ if it is not.  Second, we let $c$ be the maximal integer for which
$$(hc+1) \cdot p^k + \delta \leq hm-h+1.$$ Note that $k$ and $c$ are nonnegative integers and $c \leq p-1,$ since for $c \geq p$ we would have
$$(hc+1) \cdot p^k \geq p^{k+1} +\delta > hm-h+1.$$  It is also worth noting that $$f_1(m,h)=hm-h+1.$$

Our sufficient condition can now be stated as follows:

\begin{thm}  \label{sufficient for =}

Suppose that $2 \leq h \leq p-1$, and let $k$ and $c$ be the unique nonnegative integers defined above.  If 
$$m \leq (c+1) \cdot p^k,$$then  
$$\rho_{\pm} (\mathbb{Z}_p^r, m, h) = \rho (\mathbb{Z}_p^r, m, h).$$

\end{thm}

In fact, we believe that this condition is also necessary:

\begin{conj} \label{conj p square}

The converse of Theorem \ref{sufficient for =} is true as well; that is, if $2 \leq h \leq p-1$, $k$ and $c$ are the unique nonnegative integers defined above, and 
$$m > (c+1) \cdot p^k,$$then  
$$\rho_{\pm} (\mathbb{Z}_p^r, m, h) > \rho (\mathbb{Z}_p^r, m, h).$$
\end{conj}

We are able to prove that Conjecture \ref{conj p square} holds in the case of $\rho_{\pm} (\mathbb{Z}_p^2, m, 2)$:

\begin{thm} \label{thm p square}
Let $p$ be an odd prime and $m \leq p^2$ be a positive integer.  Then  
$$\rho_{\pm} (\mathbb{Z}_p^2, m, 2) = \rho (\mathbb{Z}_p^2, m, 2),$$ if, and only if, one of the following holds:
\begin{itemize}
  \item $m \leq p$,
  \item $m \geq (p^2+1)/2$, or
  \item there is a positive integer $c \leq (p-1)/2$ for which
  $$c \cdot p+(p+1)/2 \leq m \leq (c+1) \cdot p.$$
\end{itemize}

\end{thm}
Our proof of Theorem \ref{thm p square} involves some of the deeper methods of additive combinatorics, including Vosper's Theorem and (Lev's improvement of) Kemperman's results on critical pairs.

According to Theorem \ref{thm p square}, for a given $p$, there are exactly $(p-1)^2/4$ values of $m$ for which $\rho_{\pm} (\mathbb{Z}_p^2, m, 2)$ and $\rho(\mathbb{Z}_p^2, m, 2)$ disagree---fewer than $1/4$ of all possible values.  We have not been able to find any groups where this proportion is higher than $1/4$.

\section{The proofs of Theorems \ref{p leq h} and \ref{sufficient for =}}

In this section we establish the two sufficient conditions for the equality $$\rho_{\pm} (\mathbb{Z}_p^r, m, h) = \rho (\mathbb{Z}_p^r, m, h)$$ that we stated in Theorems \ref{p leq h} and \ref{sufficient for =}.  In order to do so, we first classify all cases with $$u_{\pm}(\mathbb{Z}_p^r,m,h)=u(p^r,m,h).$$

Let $p$ be an odd positive prime, $r \geq 2$ an integer, and $m \leq p^r$ a positive integer.  Via the (unique) base $p$ representation of $m-1$, we write $m$ as $$m=q_{r-1}p^{r-1}+\cdots+q_1  p +q_0+1,$$ where $q_{r-1}, \dots, q_0$ are all integers between $0$ and $p-1$, inclusive.  We will also need to identify three special indices:
\begin{itemize}
  \item $i_1$ denotes the largest index $i$ for which $q_i \geq 1$; if there is no such index (that is, if $m=1$), then we let $i_1=-1$.
  \item $i_2$ denotes the largest index $i$ for which $q_i \geq p/h$; if there is no such index, we let $i_2=-1$. 
  \item $i_3$ denotes the largest index $i$ for which $i>i_2$ and $$q_i=q_{i-1}=\cdots=q_{i_2+1}=(p-1)/h;$$ if there is no such index, we let $i_3=i_2$.
\end{itemize}   
We have $i_1=\lceil \log_p m \rceil -1$, and $$r-1 \geq i_1 \geq i_3 \geq i_2 \geq -1.$$

Recall that for positive integers $n$, $m$, and $h$, $$u(n,m,h)=\min\{f_d (m,h) \mid d \in D(n)\},$$
where $$f_d (m,h) = \left( h \lceil m/d \rceil -h+1 \right) \cdot d.$$  Our next proposition exhibits all values of $d$ for which $f_d(m,h)$ equals $u(p^r,m,h)$.

\begin{prop} \label{lemma f_d p}

With our notations as above, for a nonnegative integer $i$ we have $$u(p^r,m,h)=f_{p^i}(m,h)$$ if, and only if, $$i_2+1 \leq i \leq i_3+1.$$

\end{prop}

{\em Remark:}  The fact that $$u(p^r,m,h)=f_{p^{i_3+1}}(m,h)$$ was established for $h=2$ by Eliahou and Kervaire in \cite{EliKer:2001a}.

{\em Proof:}  
Given the representation of $m$ as above, we find that for every $0 \leq i \leq r$,
$$f_{p^i}=f_{p^i} (m,h)= 
h \cdot \left( q_{r-1}p^{r-1}+ \cdots + q_i p^i      \right) + p^i .$$
Therefore, when $i \geq i_3+2$, then
\begin{eqnarray*}
f_{p^i}- f_{p^{i_3+1}} & = & 
p^i-p^{i_3+1}-h \cdot \left( q_{i-1}p^{i-1}+ \cdots + q_{i_3+1} p^{i_3+1}      \right) \\
& > & p^i-p^{i_3+1}-(p-1) \cdot \left( p^{i-1}+ \cdots +  p^{i_3+1}      \right) \\
& = & 0.
\end{eqnarray*} 
Similarly, when $i_2+1 \leq i \leq i_3+1$, then
\begin{eqnarray*}
f_{p^{i_3+1}}- f_{p^i} & = & 
p^{i_3+1}-p^{i}-h \cdot \left( q_{i_3}p^{i_3}+ \cdots + q_{i} p^{i}      \right) \\
& = & p^{i_3+1}-p^{i}-(p-1) \cdot \left( p^{i-1}+ \cdots +  p^{i_3+1}      \right) \\
& = & 0.
\end{eqnarray*}
Finally, if $0 \leq i \leq i_2$, then $h q_{i_2} \geq p$, so we have
\begin{eqnarray*}
f_{p^{i}}- f_{p^{i_2+1}} & = & 
p^i-p^{i_2+1}+h \cdot \left( q_{i_2}p^{i_2}+ \cdots + q_{i} p^{i}      \right) \\
& \geq & p^i-p^{i_2+1}+h q_{i_2}p^{i_2} \\
& > & 0,
\end{eqnarray*}
completing our proof. $\Box$

Recall that for a group of type $(n_1, \dots, n_r)$ and positive integers $m$, and $h$, we defined $$u_{\pm} (G, m,h) = \min \{f_d (m,h) \mid d \in D(G,m) \},$$
where $$D(G,m)=\{d \in D(n) \mid d= d_1 \cdots d_r, d_1 \in D(n_1), \dots, d_r \in D(n_r), dn_r \geq d_rm \}.$$
Our next result finds all values of $i$ for which $f_{p^i}(m,h)$ equals $u_{\pm} (\mathbb{Z}_p^r, m,h)$.

\begin{prop} \label{u+- for p}

Let $m \geq 2$.  With our notations as above, for a nonnegative integer $i$ we have $$u_{\pm}(\mathbb{Z}_p^r, m,h)=f_{p^i}(m,h)$$ if, and only if, 
$$i=\left\{
\begin{array}{cll}
i_1+1 & \mbox{if } & hq_{i_1} \geq p; \\
i_1 \; \mbox{or } \; i_1+1 & \mbox{if } & hq_{i_1} = p-1; \\
i_1 & \mbox{if } & hq_{i_1} \leq p-2. \\
\end{array}
\right.$$

\end{prop}

{\em Proof:}  By Theorem \ref{u pm with f}, $$u_{\pm} (\mathbb{Z}_p^r, m,h) =\min \{ f_d(m,h) \mid d \in D(\mathbb{Z}_p^r,m)\}.$$ By our definition of $i_1$,     
$$D(\mathbb{Z}_p^r,m)=\{p^i \mid i_1 \leq i \leq r\}.$$  The explicit result now follows via the same considerations as in the proof of Proposition \ref{lemma f_d p} above---we omit the details.  $\Box$

Propositions \ref{lemma f_d p} and \ref{u+- for p} then imply the following:

\begin{prop}  \label{when u equals u pm}
With our notations as above, we have $$u(p^r,m,h)=u_{\pm}(\mathbb{Z}_p^r,m,h)$$ if, and only if, $i_1=i_3$ or $i_1=i_3+1$.

\end{prop}

We are now ready for the proofs of Theorems \ref{p leq h} and \ref{sufficient for =}.  

{\em Proof of Theorem \ref{p leq h}:}  Note that when $p \leq h$, then, for each index $i$, $q_i \geq 1$ is equivalent to $q_i \geq p/h$.  Therefore, we have $$i_1=i_2=i_3,$$ so our result follows from Proposition \ref{when u equals u pm} and Theorem \ref{u pm with f}.  $\Box$

{\em Proof of Theorem \ref{sufficient for =}:}  By assumption, we have nonnegative integers $k$ and $c$ with $ c \leq p-1$ so that 
$$c \cdot p^k + \frac{1}{h} \cdot (p^k-1+\delta)+1 \leq m \leq (c+1) \cdot p^k,$$ where $\delta=0$ if $p-1$ is divisible by $h$, and $\delta=1$ if it is not.  Therefore, $$i_1=\left\{
\begin{array}{cll}
k-1 & \mbox{if } & c=0, \\
k & \mbox{if } & c \geq 1. \\
\end{array}
\right.$$

To find $i_3$, assume first that $p-1$ is divisible by $h$.  Our bounds for $m$ above can then be written as 
$$c \cdot p^k +  \frac{p-1}{h} \cdot p^{k-1} +  \frac{p-1}{h} \cdot p^{k-2} + \cdots + \frac{p-1}{h} +1 \leq m \leq (c+1) \cdot p^k.$$ Thus we see that, no matter what $i_2$ equals, we have
$$i_3=\left\{
\begin{array}{cll}
k-1 & \mbox{if } & c < (p-1)/h, \\
k & \mbox{if } & c \geq (p-1)/h. \\
\end{array}
\right.$$
Our result now follows from Proposition \ref{when u equals u pm} and Theorem \ref{u pm with f}.

The case when $p-1$ is not divisible by $h$ is similar; this time the bounds for $m$ are
 $$c \cdot p^k + \frac{p}{h} \cdot p^{k-1}+1 \leq m \leq (c+1) \cdot p^k,$$ so
$$i_2=i_3=\left\{
\begin{array}{cll}
k-1 & \mbox{if } & c < p/h, \\
k & \mbox{if } & c \geq p/h; \\
\end{array}
\right.$$
implying our claim as before.  Our proof is thus complete.  $\Box$

\section{The proof of Theorem \ref{thm p square}}

We now turn to the question of determining all values of $m$ for which $$\rho_{\pm} (\mathbb{Z}_p^2, m,2) = \rho (\mathbb{Z}_p^2, m,2).$$  In order to do so, we will need to discuss some results on the so-called inverse problem in additive combinatorics; in particular, we will review some of what is known about subsets $A$ and $B$ of $G$ when their sumset $A+B$ is small.  

Recall that a subset $A$ of an abelian group $G$ is called an arithmetic progression if it is of the form $$A=\{a +i \cdot b \; : \; 0 \leq i \leq m-1\}$$ for some elements $a, b \in G$ and $m \in \mathbb{N}$; $b$ must have order at least $m$ in $G$.  Here $m$ is called the length of the arithmetic progression (we allow length 1), and $b$ is called the common difference of the progression. 

The first nontrivial inverse theorem is Vosper's classic result for groups of prime order:

\begin{thm} [Vosper; cf.~\cite{Vos:1956a} and \cite{Vos:1956b}] 
Suppose that $A$ and $B$ are nonempty subsets of $\mathbb{Z}_p$ satisfying
$$|A+B|=|A|+|B|-1.$$  
Then at least one of the following holds:
\begin{itemize}
\item $|A|=1$ or $|B|=1$;
\item $|A|+|B| = p+1$;
\item $A=\mathbb{Z}_p \setminus (g -B)$ where $\{g\} = \mathbb{Z}_p \setminus (A+B)$; or
\item $A$ and $B$ are both arithmetic progressions with the same common difference.
\end{itemize}
\end{thm}

For our use below, the following immediate consequence of Vosper's Theorem is sufficient:

\begin{cor} \label{cor to vosper}
Suppose that $A$ and $B$ are nonempty subsets of $\mathbb{Z}_p$ satisfying $|B| \geq 2$ and
$$|A+B|=|A|+|B|-1 \leq p-2.$$ 
Then $A$ is an arithmetic progression.
\end{cor}

For groups of composite order, the situation is considerably more complicated due to the existence of nontrivial proper subgroups.   Nevertheless, Kemperman \cite{Kem:1960a} gave a complete characterization of all {\em critical pairs} of finite subsets of an abelian group; that is, all finite subsets $A$ and $B$ for which $$|A+B| \leq |A|+|B|-1.$$  Kemperman's characterization was rather complicated, but it facilitated several improvements, of which we find Lev's following result most helpful for our purposes:

\begin{thm} [Lev; cf.~\cite{Lev:2006a} Theorem 4] \label{Lev's thm}

Let $A$ and $B$ be nonempty finite subsets of an abelian group $G$ satisfying $|B| \geq 2$ and
$$|A+B| \leq \min\{|G|-2,|A|+|B|-1\}.$$  Then at least one of the following holds:
\begin{itemize}
  
\item $A$ is an arithmetic progression;

 \item there exists a nonzero subgroup $H$ of $G$ with finite index $t \geq 2$ so that $A$ is the disjoint union of an arithmetic progression and $t-1$ cosets of $H$; or

\item there exists a finite, nonzero subgroup $H$ of $G$ such that 
$$|A+H| \leq \min\{|G|-1,|A|+|H|-1\}.$$

\end{itemize}

\end{thm}

We are now ready to embark on our proof of Theorem \ref{thm p square}.

{\em Proof of Theorem \ref{thm p square}:}  First, we show how the ``if'' direction follows from Theorem \ref{sufficient for =}: Keeping the notations introduced there, we see that
$$(k,c) = \left\{
\begin{array}{cll}
(0,m-1) & \mbox{if} & 1 \leq m \leq (p-1)/2; \\ \\
(1,0) & \mbox{if} & (p+1)/2 \leq m \leq p; \\ \\
(1,q) & \mbox{if} & m=qp+v \; \mbox{with } 1 \leq q \leq (p-1)/2 \; \; \mbox{and } (p+1)/2 \leq v \leq p; \\ \\
(2,0) & \mbox{if} & (p^2+1)/2 \leq m \leq p^2.
\end{array}
\right.$$
In each case we find that $m \leq (c+1) \cdot p^k.$

This leaves us with subset sizes of the form
$$m=qp+v$$ with $$1 \leq q \leq (p-1)/2 \; \; \; \mbox{and} \; \; 1 \leq v \leq (p-1)/2;$$
we will prove that, in this case,  $$\rho_{\pm} (\mathbb{Z}_p^2, m,2) > \rho(\mathbb{Z}_p^2, m,2).$$ 
Since Proposition \ref{lemma f_d p} yields 
$$\rho(\mathbb{Z}_p^2, m,2)=f_1(m,2)=2m-1,$$ our goal is to prove that
$$\rho_{\pm} (\mathbb{Z}_p^2, m,2) \geq 2m.$$

Let $A$ be an $m$-subset of $\mathbb{Z}_p^2$ for which $$|2_{\pm} A |=\rho_{\pm} (\mathbb{Z}_p^2, m,2);$$ furthermore, by Theorem \ref{symmetry thm}, we may also assume that $A$ is symmetric, near-symmetric, or asymmetric.  We will prove that $|2_{\pm} A | \geq 2m$.  

First, let us deduce what Theorem \ref{Lev's thm} says about our situation.  Following an indirect approach, let us assume that $2_{\pm} A$, and thus $2A$, have size at most $2m-1$.  (They will then have size $2m-1$.)  Note that $$2m-1 =2qp+2v-1 \leq p^2-2,$$ so the conditions of the theorem are met with $B=A$.  

Note also that $q,v \geq 1$, so $m >p$, and thus $A$ cannot be an arithmetic progression in $\mathbb{Z}_p^2$.  

Furthermore, a nonzero subgroup $H$ of index at least 2 must be of order $p$ and index $p$; with $$|A|=qp+v \leq (p^2-1)/2<(p-1)p,$$ $A$ cannot contain the disjoint union of $p-1$ distinct cosets of $H$.  

This leaves only one possibility: there must exist a subgroup $H$ of $\mathbb{Z}_p^2$ of order $p$ for which
$$|A+H| \leq \min\{p^2-1,m+p-1\}=m+p-1.$$

Now $A \subseteq A+H$, so 
$$qp < m \leq  |A+H| \leq m+p-1 =(q+1)p+v-1 < (q+2)p,$$ and thus $A+H$ is the union of exactly $q+1$ distinct cosets of $H$.  Let $A_1, \dots, A_{q+1}$ be the intersections of these cosets with $A$; $A$ is then the union of these $q+1$ components.    
 
Let us see what we can say about the sizes of these components.  Since
$$(q-1) \cdot p +2 \cdot (p-1)/2 <m,$$ at most one of the components has size less than $(p+1)/2$.  Furthermore,
$$(q-1) \cdot p+p <m$$ as well, so any two distinct components have a combined size of at least $p+1$.  Therefore, we have two possibilities:

(i) each component has size at least $(p+1)/2$; or

(ii) one component, say $A_1$, has size at most $(p-1)/2$, all other components have size at least $(p+1)/2$, and $|A_1| + |A_{i}| \geq p+1$ for all $i=2,3,\dots,q+1$.  

By the Cauchy--Davenport Theorem, applied to the $q+1$ cosets in the group $$\mathbb{Z}_p^2/H \cong \mathbb{Z}_p,$$ we can conclude that $2A$ lies in---that is, intersects non-trivially---at least $$\min\{p,2(q+1)-1\}=2q+1$$ cosets of $H$.  

Citing the Cauchy--Davenport Theorem again, we see that in case (i), $2A$ is actually the union of these cosets, and so we have
$$|2A| \geq (2q+1)p>2m,$$ contradicting our indirect assumption.

Observe that even in case (ii), with the possible exception of the coset containing $2A_1$, all the cosets that $2A$ lies in are entirely contained in $2A$.  Therefore, if $2A$ lies in at least $2q+2$ cosets of $H$, then we still have $$|2A| \geq (2q+1)p>2m,$$ contradicting our indirect assumption.

Suppose then that we are in case (ii) and that $2A$ lies in exactly $2q+1$ cosets of $H$.  Now if $|A_1| \geq v+1$, then 
$$|2A| \geq 2qp+ \min\{p,2(v+1)-1\}=2qp+2v+1>2m,$$ a contradiction again.  On the other hand, if $|A_1| \leq v$, then, since $$qp+v=|A|=|A_1|+ \cdots + |A_{q+1}|,$$ we must have $|A_1|=v$ and $|A_i| =p$ for each $i=2,3,\dots,q+1$.  Thus we are in the situation where there are exactly $q+1$ cosets intersecting $A$, one component---namely, $A_1$---has size $v$, and the other components all have size $p$.  Furthermore, $2A$ lies in exactly $2q+1$ cosets, and $2q$ of these cosets---all but the one containing $2A_1$---are entirely in $2A$.

Suppose now that $q \leq (p-3)/2$.  Then $2q+1 \leq p-2$, so by Corollary \ref{cor to vosper} of Vosper's Theorem, the $q+1$ cosets that $A$ lies in must form an arithmetic progression.  
Therefore, we can write $A$ in the form
$$A=\cup_{i=0}^{q} (a+ig+H_i),$$ where $a$ and $g$ are group elements, $H_i$ is a subset of $H$ for each $i$, and $H_i=H$ for all but one $i$.   Consequently, for distinct $i_1$ and $i_2$, at least one of $H_{i_1}$ or $H_{i_2}$ equals $H$, so we have
$$(a+i_1g+H_{i_1})-(a+i_2g+H_{i_2})=(i_1-i_2)g+H \subseteq 2_{\pm}A.$$  Furthermore, since $q +1 \geq 2$, there is an $i$ such that $H_i=H$, so for this $i$ we have $$(a+ig+H_i)-(a+ig+H_i)=H,$$ and thus $$H \setminus \{0\} \subseteq 2_{\pm}A.$$  This implies that
$$\cup_{i=-q}^{q} (ig+H) \setminus \{0\} \subseteq 2_{\pm}A,$$ and so
$$|2_{\pm}A| \geq (2q+1)p-1 \geq 2m,$$ a contradiction with our indirect assumption.  

We are left with the case when $q=(p-1)/2$, in which case
\begin{itemize}
  \item $A$ has size $(p^2-p)/2+v$ and is the union of $A_1$ and $(p-1)/2$ cosets of $H$, and
\item $2A$ has size $p^2-p+2v-1$ and is the union of $2A_1$ and $p-1$ cosets of $H$.
\end{itemize}  
Recall that we are assuming that $A$ is symmetric, near-symmetric, or asymmetric---we attend to each of these cases separately.

Suppose first that $A$ is near-symmetric.  Then, by definition, we can find an element $a \in \mathbb{Z}_p^2 \setminus A$ for which $$A'=A \cup \{a\}$$ is symmetric.  We can easily check that we then have $$2_{\pm} A'=2_{\pm}A$$
(note that $0 \in 2_{\pm}A$ since $m \geq 2$).  Therefore, 
$$2m-1 \geq |2_{\pm}A| = |2_{\pm}A'| \geq \rho_{\pm} (\mathbb{Z}_p^2,m+1,2) \geq \rho (\mathbb{Z}_p^2,m+1,2).$$
However, this is a contradiction, since by Proposition \ref{lemma f_d p}, we get
$$\rho (\mathbb{Z}_p^2,m+1,2)=\left\{
\begin{array}{ll}
f_1(m+1,2)=2(m+1)-1=2m+1 & \mbox{if} \; v \leq (p-3)/2; \\ \\
f_p(m+1,2)=p^2 = 2m+1 & \mbox{if} \; v = (p-1)/2.
\end{array}\right.$$

Next, assume that $A$ is symmetric.  Observe that we then must have $A_1 \subset H$, since otherwise $A$ would contain fewer than $(p-1)/2$ full cosets of $H$.  

Now suppose that $a+H$ is one of the cosets that make up $A$.  Then, because it is symmetric, $A$ also contains $-a+H$; and since these are two disjoint cosets, $2A$ will contain their sum, which is $H$.  But this then implies that $2A=\mathbb{Z}_p^2$, which again contradicts $|2A| \leq 2m-1$.

Finally, suppose that $A$ is asymmetric.  Consider the partition 
$$H \cup (\pm a_1 +H) \cup \cdots \cup (\pm a_{(p-1)/2} +H)$$ of $\mathbb{Z}_p^2$ (here $a_1, \dots, a_{(p-1)/2}$ are appropriately chosen group elements).  Recall that $(p-1)/2$ of these $p$ cosets must lie entirely in $A$, so the fact that $A$ is asymmetric implies that $A_1 \subset H$.  With $a+H$ being one of the cosets in $A$,
$$(a+H) - (a+H) =H,$$ and thus $$H \setminus \{0\} \subseteq 2_{\pm}A.$$  Therefore,   
$$|2_{\pm}A| \geq (p-1)p+(p-1) =p^2-1 \geq 2m,$$ a contradiction.

Our proof is now complete. $\Box$

We have thus identified each value of $m$ for which $\rho_{\pm} (\mathbb{Z}_p^2,m,2)$ equals $\rho (\mathbb{Z}_p^2,m,2)$---and thus also equals $u(p^2,m,2)$---and those for which it does not, but how about the exact value of $\rho_{\pm} (\mathbb{Z}_p^2,m,2)$?  We believe that  
$\rho_{\pm} (\mathbb{Z}_p^2,m,2)$ equals $u_{\pm}(\mathbb{Z}_p^2,m,2)$ for all values of $m$, except for when
 $$\frac{p^2-p+2}{2} \leq m \leq \frac{p^2-1}{2}.$$  Indeed, in this case we have
$$u_{\pm}(\mathbb{Z}_p^2,m,2)=p^2$$ (cf. Proposition \ref{u+- for p}), but
$$\rho_{\pm} (\mathbb{Z}_p^2,m,2) \leq p^2-1,$$ as demonstrated by any asymmetric $m$-subset $A$ of $\mathbb{Z}_p^2$ (we then have $0 \not \in 2_{\pm} A$).

In light of this, we make the following conjecture: 
\begin{conj}

Let us write $m$ as $m=cp+v$ with $$0 \leq c \leq p-1 \; \; \; \mbox{and} \; \; \; 1 \leq v \leq p.$$  We then have:
$$\begin{array}{cc|ccccc} 
c & v & \rho (\mathbb{Z}_p^2,m,2) && \rho_{\pm} (\mathbb{Z}_p^2,m,2) && u_{\pm}(\mathbb{Z}_p^2,m,2) \\ \hline 
& & & & \\
 & v \leq (p-1)/2 & 2m-1 &=& 2m-1 &=& 2m-1 \\ 
\raisebox{1.5ex}[0pt]{$0$} & v \geq (p+1)/2 & p & =&p & =&p \\ 
& & & & \\
 & v \leq (p-1)/2 & 2m-1 & <& \fbox{$(2c+1)p$}  & =&(2c+1)p \\ 
\raisebox{1.5ex}[0pt]{$1 \leq c \leq (p-3)/2$} & v \geq (p+1)/2 & (2c+1)p & =&(2c+1)p & =&(2c+1)p \\  
& & & & \\
 & v \leq (p-1)/2 & 2m-1 & <& \fbox{$p^2-1 $} & <& p^2 \\ 
\raisebox{1.5ex}[0pt]{$c=(p-1)/2$} & v \geq (p+1)/2 & p^2 & = & p^2 & =&p^2 \\ 
& & & & \\
c \geq (p+1)/2 & \mbox{any } v & p^2 & =&p^2 & =&p^2 \\ 
\end{array}$$

\end{conj}

The two boxed entries in this table remain unproven.


\begin{thebibliography}{99}


\bibitem{Baj:2000a} B. Bajnok, Spherical Designs and Generalized Sum-Free Sets in Abelian Groups.  Special issue dedicated to Dr. Jaap Seidel on the occasion of his 80th birthday (Oisterwijk, 1999).  \emph{Des. Codes Cryptogr.}  {\bf 21}(1-3) (2000), 11-18.

\bibitem{Baj:2004a} B. Bajnok, The Spanning Number and the Independence Number of a Subset of an Abelian Group.  In \emph{Number Theory,} D. Chudnovsky, G. Chudnovsky, and M. Nathalson (Ed.), Springer-Verlag (2004), 1-16.

\bibitem{BajMat:2014a} B. Bajnok and R. Matzke, The Minimum Size of Signed Sumsets, www.arxiv.org (2014).



\bibitem{BajRuz:2003a} B. Bajnok and I. Ruzsa, The Independence Number of a Subset of an Abelian Group.  \emph{Integers} {\bf 3}(A2) (2003), 23 pp.  

\bibitem{Cau:1813a} A.-L. Cauchy, Recherches sur les nombres, {\em J. \'Ecole Polytechnique} {\bf 9} (1813) 99--123.


\bibitem{Dav:1935a} H. Davenport, On the addition of residue classes, {\em J. London Math. Soc.} {\bf 10} (1935) 30--32.

\bibitem{Dav:1947a} H. Davenport, A historical note, {\em J. London Math. Soc.} {\bf 22} (1947) 100--101.

\bibitem{EliKer:2001a} S. Eliahou and M. Kervaire, Restricted Sumsets in Finite Vector Spaces: The Case $p=3$, {\em Integers}, {\bf 1} (2001) \#A02.

\bibitem{EliKer:2005a} S. Eliahou and M. Kervaire, Old and new formulas for the Hopf--Stiefel and related functions, {\em Expo. Math.}, {\bf 23}:2 (2005) 127--145.

\bibitem{EliKer:2007a} S. Eliahou and M. Kervaire, Some extensions of the Cauchy--Davenport Theorem, {\em Electronic Notes in Discrete Math.}, {\bf 28} (2007) 557--564.

\bibitem{EliKerPla:2003a} S. Eliahou, M. Kervaire, and A. Plagne, Optimally small sumsets in finite abelian groups, {\em J. Number Theory}, {\bf 101} (2003) 338--348.


\bibitem{Kar:2006a} Gy. K\'arolyi, A note on the Hopf--Stiefel function.  {\em European J. Combin.}, {\bf 27} (2006) 1135--1137.

\bibitem{Kem:1960a} J. H. B. Kemperman, On small sumsets in an abelian group, {\em Acta Mathematica}, {\bf 103} (1960) 63--88.

\bibitem{KloLev:2003a} B. Klopsch and V. F. Lev, How long does it take to generate a group?  {\em J. Algebra}, {\bf 261}, (2003) 145--171.

\bibitem{KloLev:2009a} B. Klopsch and V. F. Lev, Generating abelian groups by addition only.  {\em Forum Math.}, {\bf 21}:1, (2009) 23--41.

\bibitem{Lev:2006a} V. F. Lev, Critical pairs in abelian groups and Kemperman's structure theorem, {\em Int. J. Number Theory} {\bf 3} (2006) 379--396.

\bibitem{Pla:2003a} A. Plagne, Additive number theory sheds extra light on the Hopf--Stiefel $\circ$ function, {\em Enseign. Math., II S\'er}, {\bf 49}:1--2 (2003) 109--116.

\bibitem{Pla:2006a} A. Plagne, Optimally small sumsets in groups, I. The supersmall sumset property, the $\mu_G^{(k)}$ and the $\nu_G^{(k)}$ functions, {\em Unif. Distrib. Theory}, {\bf 1}:1 (2006) 27--44.

\bibitem{Sha:1984a} D. Shapiro, Products of sums of squares, {\em Expo. Math.}, {\bf 2} (1984) 235--261.




\bibitem{Vos:1956a} A. G. Vosper, The critical pairs of subsets of a group of prime order.  {\em J. London Math. Soc.} {\bf 31} (1956) 200--205; 

\bibitem{Vos:1956b} A. G. Vosper, Addendum to ``The critical pairs of subsets of a group of prime order.''  {\em J. London Math. Soc.} {\bf 31} (1956) 280--282.






\end{thebibliography}
\end{document}